\documentclass[12pt,leqno,a4paper,twoside]{amsart}
\usepackage{amsmath,amsthm,amssymb,verbatim,enumerate,ifthen,times,cite}
\usepackage[mathscr]{eucal}
\usepackage[utf8]{inputenc}
\usepackage[T1]{fontenc}

\newcommand{\N}{\mathbb{N}}
\newcommand{\Z}{\mathbb{Z}}
\newcommand{\R}{\mathbb{R}}
\newcommand{\A}{\mathscr{A}}
\newcommand{\C}{\mathscr{C}}
\newcommand{\T}{\mathscr{T}}
\renewcommand{\S}{\mathscr{S}}

\newtheorem{theorem}{Theorem}[section]
\newtheorem*{theorem*}{Theorem}
\def\Thm#1#2{\ifthenelse{\equal{#1}{*}}{\begin{theorem*}#2\end{theorem*}}
  {\begin{theorem}\label{T#1}#2\end{theorem}}}
\newtheorem{Atheorem}{Theorem}

\def\thm#1{Theorem~\ref{T#1}}

\newtheorem{proposition}[theorem]{Proposition}
\newtheorem*{proposition*}{Proposition}
\def\Prp#1#2{\ifthenelse{\equal{#1}{*}}{\begin{proposition*}#2\end{proposition*}}
             {\begin{proposition}\label{P#1}#2\end{proposition}}}
\def\prp#1{Proposition~\ref{P#1}}

\newtheorem{corollary}[theorem]{Corollary}
\newtheorem*{corollary*}{Corollary}
\def\Cor#1#2{\ifthenelse{\equal{#1}{*}}{\begin{corollary*}#2\end{corollary*}}
             {\begin{corollary}\label{C#1}#2\end{corollary}}}

\newtheorem{lemma}[theorem]{Lemma}
\newtheorem*{lemma*}{Lemma}
\def\Lem#1#2{\ifthenelse{\equal{#1}{*}}{\begin{lemma*}#2\end{lemma*}}
             {\begin{lemma}\label{L#1}#2\end{lemma}}}
\def\lem#1{Lemma~\ref{L#1}}

\newtheorem{example}{Example}
\newtheorem*{example*}{Example}
\def\Exa#1#2{\ifthenelse{\equal{#1}{*}}{\begin{example*}\rm #2\end{example*}}
             {\begin{example}\label{Ex#1}\rm #2\end{example}}}

\newtheorem{problem}[theorem]{Problem}

\theoremstyle{definition}
\newtheorem{definition}[theorem]{Definition}
\newtheorem*{definition*}{Definition}
\def\Defi#1#2{\ifthenelse{\equal{#1}{*}}{\begin{definition*}#2\end{definition*}}
             {\begin{definition}\label{D#1}#2\end{definition}}}

\newtheorem{remark}[theorem]{Remark}
\newtheorem*{remark*}{Remark}
\def\Rem#1#2{\ifthenelse{\equal{#1}{*}}{\begin{remark*}\rm #2\end{remark*}}
             {\begin{remark}\label{R#1}\rm #2\end{remark}}}

\newcommand{\eq}[1]{\eqref{E#1}}
\def\Eq#1#2{\ifthenelse{\equal{#1}{*}}
  {\begin{equation*}\begin{aligned}[]#2\end{aligned}\end{equation*}}
  {\begin{equation}\begin{aligned}[]\label{E#1}#2\end{aligned}\end{equation}}}

\begin{document}

\title{On convex and concave sequences and their applications}

\author[G. M. Molnár]{Gábor M. Molnár}
\address{Doctoral School of Mathematical and Computational Sciences, University of Debrecen, H-4002 Debrecen, Pf.\ 400, Hungary}
\email{molnar.gabor.marcell@science.unideb.hu}

\author[Zs. P\'ales]{Zsolt P\'ales}
\address{Institute of Mathematics, University of Debrecen, H-4002 Debrecen, Pf.\ 400, Hungary}
\email{pales@science.unideb.hu}

\subjclass[2010]{26A51, 39B62}
\keywords{$q$-convex sequence, $q$-concave sequence, $q$-affine sequence, Chebyshev polynomials of the first and second kind, contraction}


\thanks{The research of the second author was supported by the K-134191 NKFIH Grant and the 2019-2.1.11-T\'ET-2019-00049, EFOP-3.6.1-16-2016-00022 and EFOP-3.6.2-16-2017-00015 projects. The last two projects are co-financed by the European Union and the European Social Fund.}

\begin{abstract}
The aim of this paper is to introduce and to investigate the basic properties of $q$-convex, $q$-affine and $q$-concave sequences and to establish their surprising connection to Chebyshev polynomials of the first and of the second kind. One of the main results shows that $q$-concave sequences are the pointwise minima of $q$-affine sequences. As an application, we consider a nonlinear selfmap of the $n$-dimensional space and prove that it has a unique fixed point. For the proof of this result, we introduce a new norm on the space in terms of a $q$-concave sequence and show that the nonlinear operator becomes a contraction with respect to this norm, and hence, the Banach Fixed Point theorem can be applied.
\end{abstract}

\maketitle

\section{Introduction}

In the theory of convexity, the investigation of convex functions play a fundamental role. We refer to the following monographs for the details: Hardy--Littlewood--P\'olya \cite{HarLitPol34}, Kuczma \cite{Kuc85}, Mitrinovi\'c \cite{Mit70}, Mitrinovi\'c--Pe\v{c}ar\'c--Fink \cite{MitPecFin91,MitPecFin93}, Niculescu--Persson \cite{NicPer06}, Popoviciu \cite{Pop44}, and Roberts--Varberg \cite{RobVar73}. The investigation of convex sequences probably started in the book Mitrinovi\'c \cite{Mit70}. This subfield is still very active, some recent results and applications have been obtained by Krasniqi \cite{Kra16}, Niezgoda \cite{Nie11,Nie17a,Nie17b}, Sofonoea--\c{T}incu--Acu \cite{SofTinAcu18}, Tabor--Tabor--\.Zoldak \cite{TabTabZol12}, Wu--Debnath \cite{WuDeb07}, Y\i ld\i z \cite{Yil18}. In this paper we introduce the notions of $q$-convex, $q$-affine and $q$-concave sequences and we present some basic results on them and we establish their surprising connection to Chebyshev polynomials of the first and of the second kind. Finally, we present an application of them to fixed point theory.

Let $\R$, $\R_+$, $\Z$ and $\N$ denote the sets of real, positive real, integer and positive integer numbers in this paper. Given $n,m\in\Z$ with $2\leq m-n$, let $\S(n|m)$ denote the linear space $\R^{\{n,\dots,m\}}$ of all real sequences, i.e., the collection of all functions $p:\{n,\dots,m\}\to\R$. It is natural to define the notions of concavity, convexity and affinity for the elements of $\S(n|m)$. A sequence $p=(p_n,\dots,p_m)\in\S(n|m)$ is called \emph{convex} if, for all $i\in\{n+1,\dots,m-1\}$,
\Eq{J}{
  2p_i\leq p_{i-1}+p_{i+1}.
}
If, for all $i\in\{n+1,\dots,m-1\}$, the reversed inequality holds in \eq{J}, then the sequence is termed \emph{concave}. Finally, if a sequence is simultaneously convex and concave, then it is said to be \emph{affine}. If the inequality \eq{J} holds with strict inequality sign, then we speak about strict convexity and concavity, respectively. 

In what follows, we extend the above definitions and introduce the notions of $q$-convex, $q$-concave, and $q$-affine sequences with respect to a positive number $q$. A sequence $p=(p_n,\dots,p_m)\in\S(n,m)$ is called \emph{$q$-convex} if, for all $i\in\{n+1,\dots,m-1\}$, 
\Eq{qJ}{
  2qp_i\leq p_{i-1}+p_{i+1}.
}
If, for all $i\in\{n+1,\dots,m-1\}$, the reversed inequality holds in \eq{qJ}, then the sequence is termed \emph{$q$-concave}. If a sequence is simultaneously $q$-convex and $q$-concave, then it is said to be \emph{$q$-affine}. 

We can easily see that the strict convexity of a positive (or negative) sequence implies its $q$-convexity for some $q$. Indeed, if $p\in\S(nm)$ is a positive strictly convex sequence then, for all $i\in\{n+1,\dots,m-1\}$, 
\Eq{*}{
  1<\frac{p_{i-1}+p_{i+1}}{2p_i}.
}
Therefore,
\Eq{*}{
  1<q:=\min_{i\in\{n+1,\dots,m-1\}}
  \frac{p_{i-1}+p_{i+1}}{2p_i},
}
which implies that $p$ is $q$-convex with a number $q>1$.
Analogously, $p\in\S(nm)$ is a negative strictly convex sequence, then it is $q$-convex with a number $0<q<1$.

The subclasses of $q$-convex and $q$-concave sequences in $\S(n|m)$ will be denoted $\C^\cup_q(n|m)$ and $\C^\cap_q(n|m)$, respectively. Finally, $\A_q(n|m)$ will stand for the subclass of $q$-affine sequences, that is, 
\Eq{*}{
\A_q(n|m):=\C^\cup_q(n|m)\cap\C^\cap_q(n|m).
}
It is easy to see that $\A_q(n|m)$ is a linear subspace of $\S(n|m)$ and $\C^\cup_q(n|m)$ and $\C^\cap_q(n|m)$ are convex cones in $\S(n|m)$, i.e., they are closed with respect linear combinations with nonnegative coefficients.

The aim of this paper is to investigate the basic properties of these classes of sequences and to show their surprising connection to Chebyshev polynomials of the first and of the second kind. Therefore, in the next section, we recall the notions of Chebyshev polynomials and establish the basic relationships among them.

In Section 3, we describe all $q$-affine sequences in terms of Chebyshev polynomials and show that $\A_q(n|m)$ is a two-dimensional linear subspace of $\S(n|m)$. In another result of this section, we deduce inequalities that are consequences of the $q$-convexity/concavity and we also establish an analogue of the so called support theorem and thus we obtain that $q$-concave sequences are the pointwise minima of $q$-affine sequences. 

In Section 4, we consider minimum problems for positive sequences in terms of a (power) mean $M$. In the cases when $M$ is either the arithmetic, or the geometric, or the maximum mean we obtain the precise solution of this minimum problem. For a general power mean with a positive parameter, we only obtain lower bounds. The case when $M=\max$ is strongly connected to the results obtained for $q$-concave sequences.

In the last section, we consider a nonlinear selfmap of the $n$-dimensional space $\R^n$ and prove that it has a unique fixed point. For the proof of this result, we introduce a new norm in terms of $q$-concave sequences and show that the nonlinear operator becomes a contraction with respect to this norm, and hence, by the Banach Fixed Point theorem, it has a unique fixed point.

\section{Auxiliary results for Chebyshev polynomials}

For $k\in\Z$, let $T_k:\R\to\R$ and $U_k:\R\to\R$ denote the Chebyshev polynomials of the first and of the second kind of order $k$, which are defined by the system of equations
\Eq{T}{
  T_0(x)&:=1,&\qquad T_1(x)&:=x, &\qquad 
  T_{k-1}(x)+T_{k+1}(x)&=2xT_k(x) &\qquad(k\in\Z),\\
  U_0(x)&:=1,&\qquad U_1(x)&:=2x,& \qquad 
  U_{k-1}(x)+U_{k+1}(x)&=2xU_k(x) &\qquad(k\in\Z),
}
respectively. The last equalities in \eq{T} rewritten as
\Eq{*}{
  T_{k+1}(x)=2xT_k(x)-T_{k-1}(x),\qquad 
  U_{k+1}(x)=2xU_k(x)-U_{k-1}(x),
}
can be used to compute $T_k$ and $U_k$ for $k\geq2$ recursively. If we rewrite them as
\Eq{*}{
  T_{k-1}(x)=2xT_k(x)-T_{k+1}(x),\qquad 
  U_{k-1}(x)=2xU_k(x)-U_{k+1}(x),
}
then $T_k$ and $U_k$ can be determined for $k\leq-1$.
One can easily prove that, for $k\in\Z$,
\Eq{*}{
  T_{-k}=T_k \qquad\mbox{and}\qquad U_{-k}=-U_{k-2}.
}
In particular, $U_{-1}=0$. It is clear that, for $k\geq 0$, the degree of $T_k$ and $U_k$ equals $k$. 
It is well-known that these polynomials satisfy the equalities
\Eq{TU}{
  T_k(\cos(u))=\cos(ku)\qquad\mbox{and}\qquad
   T_k(\cosh(u))=\cosh(ku) \qquad(u\in\R,\,k\in\Z)
}
and
\Eq{TU+}{
  U_{k}(\cos(u))=\frac{\sin((k+1)u)}{\sin(u)}
  \qquad\mbox{and}\qquad
  U_{k}(\cosh(u))=\frac{\sinh((k+1)u)}{\sinh(u)}
   \qquad(u\in\R,\,k\in\Z).
}
From these representations it easily follows that the roots of $T_k$ (for $k\neq0$) and $U_{k-1}$ (for $k\not\in\{-1,0,1\}$) are given by
\Eq{*}{
  \bigg\{\cos\bigg(\frac{2i-1}{2k}\pi\bigg)
  \mid i\in\{1,\ldots,|k|\}\bigg\}
  \qquad\mbox{and}\qquad
  \bigg\{\cos\bigg(\frac{i}{k}\pi\bigg)
  \mid i\in\{1,\ldots,|k|-1\}\bigg\},
}
respectively. Therefore, the largest root of $T_k$ (for $k\neq0$) and $U_{k-1}$ (for $k\not\in\{-1,0,1\}$) are given by
\Eq{*}{
  \cos\bigg(\frac{\pi}{2k}\bigg)
  \qquad\mbox{and}\qquad
  \cos\bigg(\frac{\pi}{k}\bigg),
}
respectively.

\Lem{MT}{For $0\leq x<1$, the sequence $(T_k(x))_{k=1}^{\tau(x)}$ is strictly decreasing, where $\tau(x):=\big\lfloor\frac{\pi}{\arccos(x)}\big\rfloor$. For $x>1$, the sequence $(T_k(x))_{k=0}^\infty$ is strictly increasing. }

\begin{proof} If $x>1$, then there exists $u>0$ such that $x=\cosh(u)$. Thus, in view of the second formula in \eq{TU}, we have
\Eq{*}{
  T_k(x)=T_k(\cosh(u))=\cosh(ku) \qquad(k\in\N\cup\{0\}),
}
which by the strict monotonicity of the $\cosh$ function implies that the right hand side is a strictly increasing function of $k$.

If $0\leq x<1$, then there exists $u\in\,]0,\frac{\pi}{2}]$ such that $x=\cos(u)$. In view of the first formula in \eq{TU}, we have
\Eq{*}{
  T_k(x)=T_k(\cos(u))=\cos(ku) \qquad(k\in\N\cup\{0\}),
}
which, using that $\cos$ is strictly decreasing on $[0,\pi]$, implies that $T_k(x)$ is strictly decreasing for $k\in\{0,\dots,\big\lfloor\frac{\pi}{u}\big\rfloor\}$.
\end{proof}

\Lem{Id}{Let $n\in\N$ be an odd number. Then, for all $x_1,\dots,x_n\in\R$ with the notation $x_{i+n}:=x_i$ ($i\in\{1,\dots,n\}$), we have
\Eq{Id}{
  \qquad\sum_{i=1}^n \sin\bigg(\sum_{j=1}^{n-1}(-1)^jx_{i+j}\bigg)\sin(x_i)=0
  \qquad\text{and}\qquad
  \sum_{i=1}^n \sin\bigg(\sum_{j=1}^{n-1}(-1)^jx_{i+j}\bigg)\cos(x_i)=0.
}}

\begin{proof} For $n=1$, the statement is trivial, thus we may assume that $n\geq3$. Let $x_1,\dots,x_n\in\R$ and denote
\Eq{*}{
  y_i:=\sum_{j=1}^{n-1}(-1)^jx_{i+j} \qquad(i\in\{1,\dots,n-1\}).
}
Then, by the well-known product-to-sum identities
\Eq{*}{
  2\sin\bigg(\sum_{j=1}^{n-1}(-1)^jx_{i+j}\bigg)\sin(x_i)
  &=2\sin(x_i)\sin(y_i)=\cos(x_i-y_i)-\cos(x_i+y_i),\\
  2\sin\bigg(\sum_{j=1}^{n-1}(-1)^jx_{i+j}\bigg)\cos(x_i)
  &=2\cos(x_i)\sin(y_i)=\sin(x_i+y_i)-\sin(x_i-y_i).
}
Observe that, by the equality $x_i=x_{i+n}$ and by the oddness of $n$, we have
\Eq{*}{
  x_i-y_i&=x_i-\sum_{j=1}^{n-1}(-1)^jx_{i+j}
  =x_i+x_{i+1}+\sum_{j=2}^{n-1}(-1)^{j-1}x_{i+j}\\
  &=x_{i+1}+(-1)^{n-1}x_{i+n}+\sum_{j=1}^{n-2}(-1)^jx_{i+1+j}
  =x_{i+1}+y_{i+1}.
}
Therefore 
\Eq{*}{
  2\sin\bigg(\sum_{j=1}^{n-1}(-1)^jx_{i+j}\bigg)\sin(x_i)
  &=\cos(x_{i+1}+y_{i+1})-\cos(x_i+y_i),\\
  2\sin\bigg(\sum_{j=1}^{n-1}(-1)^jx_{i+j}\bigg)\cos(x_i)
  &=\sin(x_i+y_i)-\sin(x_{i+1}+y_{i+1}).
}
Summing up these equalities side by side for $i\in\{1,\dots,n\}$, respectively, we can see that the right hand sides are telescopic sums which are equal to zero, hence both equalities in \eq{Id} hold true.
\end{proof}

\Lem{IdU}{For all $i,j,k\in\Z$, we have 
\Eq{IdU}{
  U_{k-j-1}U_i+U_{j-i-1}U_k=U_{k-i-1}U_j
  \qquad\text{and}\qquad
  U_{k-j-1}T_i+U_{j-i-1}T_k=U_{k-i-1}T_j.
}
Furthermore, for $i,j\in\Z$, we also have
\Eq{UT}{
  U_{i-j}+U_{i+j}=2T_jU_i
  \qquad\text{and}\qquad
  T_{i-j}+T_{i+j}=2T_jT_i.
}
}

\begin{proof}
In the particular case $n=3$, with $x_1:=x$, $x_2:=y$ and $x_3:=z$, the identities in \eq{Id} yield
\Eq{Id3}{
  \sin(z-y)\sin(x)+\sin(y-x)\sin(z)&=\sin(z-x)\sin(y),\\
  \sin(z-y)\cos(x)+\sin(y-x)\cos(z)&=\sin(z-x)\cos(y).
}
Let $q\in\,]-1,1[\,$ be arbitrary, let $u:=\arccos(q)$ and let $i,j,k\in\Z$. With the substitutions $(x,y,z):=((i+1)u,(j+1)u,(k+1)u)$ and $(x,y,z):=(iu,j u,ku)$, the first and second identities in \eq{Id3} imply
\Eq{*}{
  \frac{\sin((k-j)u)}{\sin(u)}\cdot\frac{\sin((i+1)u)}{\sin(u)}
  +\frac{\sin((j-i)u)}{\sin(u)}\cdot\frac{\sin((k+1)u)}{\sin(u)}
  &=\frac{\sin((k-i)u)}{\sin(u)}\cdot\frac{\sin((j+1)u)}{\sin(u)},\\
  \frac{\sin((k-j)u)}{\sin(u)}\cdot\cos(iu)
  +\frac{\sin((j-i)u)}{\sin(u)}\cdot\cos(ku)
  &=\frac{\sin((k-i)u)}{\sin(u)}\cdot\cos(ju).
}
In view of \eq{TU}, from these equalities we can easily obtain that
\Eq{*}{
  U_{k-j-1}(q)U_i(q)+U_{j-i-1}(q)U_k(q)&=U_{k-i-1}(q)U_j(q),\\
  U_{k-j-1}(q)T_i(q)+U_{j-i-1}(q)T_k(q)&=U_{k-i-1}(q)T_j(q)
}
hold for all $q\in\,]-1,1[\,$ and hence for all $q\in\R$. This completes the proof of the equalities in \eq{IdU}.

To prove \eq{UT}, let $q\in\,]-1,1[\,$ be arbitrary, let $u:=\arccos(q)$ and $i,j\in\Z$. Using \eq{TU} and the addition formula for the sine and cosine functions, we obtain
\Eq{*}{
  U_{i-j}(q)+U_{i+j}(q)&=U_{i-j}(\cos(u))+U_{i+j}(\cos(u))
  =\frac{\sin((i-j+1)u)}{\sin(u)}
  +\frac{\sin((i+j+1)u)}{\sin(u)}\\
  &=2\frac{\sin((i+1)u)}{\sin(u)}\cos(ju)
  =2U_i(\cos(u))T_j(\cos(u))=2U_i(q)T_j(q)
}
and
\Eq{*}{
  T_{i-j}(q)+T_{i+j}(q)&=T_{i-j}(\cos(u))+T_{i+j}(\cos(u))
  =\cos((i-j)u)+\cos((i+j)u)\\
  &=2\cos(iu)\cos(ju)
  =2T_i(\cos(u))T_j(\cos(u))=2T_i(q)T_j(q).
}
This completes the proof of \eq{UT}.
\end{proof}

Observe that, in the particular case $j=1$, the equalities in \eq{UT} reduce to the recursive formulas in \eq{T}

\Rem{UT+}{For the difference of two Chebyshev polynomials of the second kind, using the equality $-U_k=U_{-k-2}$, we can deduce the following identity:
\Eq{UT1}{
  U_{i+j}-U_{i-j}
  =U_{i+j}+U_{-i+j-2}
  =U_{(j-1)+(i+1)}+U_{(j-1)-(i+1)}
  =2T_{i+1}U_{j-1}.
}
On the other hand, to compute the difference of two Chebyshev polynomials of the first kind, the following equality can be established:
\Eq{UT2}{
 T_{j-i}(q)-T_{j+i}(q)=2(1-q^2)U_{j-1}(q)U_{i-1}(q).
}
To prove this, let $q\in\,]-1,1[\,$ be arbitrary, let $u:=\arccos(q)$ and $i,j\in\Z$. Using \eq{TU} and the addition formula for the cosine function, we get
\Eq{*}{
  2U_{j-1}(q)U_{i-1}(q)
  &=2U_{j-1}(\cos(u))U_{i-1}(\cos(u))
  =2\frac{\sin(ju)}{\sin(u)}\frac{\sin(iu)}{\sin(u)}\\
  &=\frac{\cos((j-i)u)-\cos((j+i)u)}{\sin^2(u)}
  =\frac{T_{j-i}(q)-T_{j+i}(q)}{1-\cos^2(u)}
  =\frac{T_{j-i}(q)-T_{j+i}(q)}{1-q^2}.
}
From here, \eq{UT2} directly follows.}

\section{$q$-concave, convex and affine sequences}

The next proposition shows that $\A_q(n|m)$ is a two dimensional subspace of $\S(n|m)$.

\Prp{qJA}{A sequence $p\in\S(n|m)$ is $q$-affine if and only if there exist $a,b\in\R$ such that
\Eq{pi}{
  p_i:=aU_{i-n}(q)+bT_{i-n}(q) \qquad(i\in\{n,\dots,m\}).
}
In addition, if $p\in\A_q(n|m)$, then, for all $i,j,k\in\{n,\dots,m\}$,
\Eq{ijk=}{
  U_{k-j-1}(q)p_{i}+U_{j-i-1}(q)p_{k}=U_{k-i-1}(q)p_j.
}
In particular, for $i\in\{n,\dots,m\}$ and $j\in\{1,\dots,\min(i-n,m-i)\}$,
\Eq{ijT}{
  p_{i-j}+p_{i+j}=2T_{j}(q)p_i.
} 
}

\begin{proof} First assume that $p=(p_n,\dots,p_m)$ is $q$-affine. Define
\Eq{*}{
  a:=\frac{p_{n+1}}{q}-p_n,\qquad b:=2p_n-\frac{p_{n+1}}{q}.
}
We prove the equality \eq{pi} by induction with respect to $i$. Observe that $p_n=a+b=aU_0(q)+bT_0(q)$ and $p_{n+1}=a(2q)+bq=aU_1(q)+bT_1(q)$, which show that \eq{pi} holds for $i=n$ and $i=n+1$. Assume that we have proved \eq{pi} for 
$i\leq\ell$, where $n+1\leq\ell\leq m-1$. Then, applying the $q$-affinity of the sequence, the inductive hypothesis and finally the recursive property of Chebyshev polynomials, we obtain
\Eq{*}{
  p_{\ell+1}=2qp_\ell-p_{\ell-1}
  &=2q(aU_{\ell-n}(q)+bT_{\ell-n}(q))-(aU_{\ell-1-n}(q)+bT_{\ell-1-n}(q))\\
  &=a(2qU_{\ell-n}(q)-U_{\ell-1-n}(q))+b(2qT_{\ell-n}(q)-T_{\ell-1-n}(q))\\
  &=aU_{\ell+1-n}(q)+bT_{\ell+1-n}(q).
}
This shows the validity of \eq{pi} for $i=\ell+1$.

For the sufficiency part of our assertion, assume that \eq{pi} holds for some $a,b\in\R$. Then, by the recursive property of Chebyshev polynomials, for $i\in\{n+1,\dots,m-1\}$, we have that
\Eq{*}{
  p_{i+1}=aU_{i+1-n}(q)+bT_{i+1-n}(q)
  &=a(2qU_{i-n}(q)-U_{i-1-n}(q))+b(2qT_{i-n}(q)-T_{i-1-n}(q))\\
  &=2q(aU_{i-n}(q)+bT_{i-n}(q))-(aU_{i-1-n}(q)+bT_{i-1-n}(q))\\
  &=2qp_i-p_{i-1},
}
which proves that $p$ is a $q$-affine sequence.

To verify the last two assertions let $p\in\A_q(n|m)$. Then, as we have seen it, \eq{pi} holds for some $a,b\in\R$. 

Let first $i,j,k\in\{n,\dots,m\}$ be arbitrary. Then, applying \lem{IdU}, we get
\Eq{*}{
  U_{k-j-1}(q)U_{i-n}(q)+U_{j-i-1}(q)U_{k-n}(q)&=U_{k-i-1}(q)U_{j-n}(q)
  \qquad\text{and}\\
  U_{k-j-1}(q)T_{i-n}(q)+U_{j-i-1}(q)T_{k-n}(q)&=U_{k-i-1}T_{j-n}(q).
}
Multiplying the first and second equalities by $a$ and $b$, respectively, and then adding them up side by side, we obtain
\Eq{*}{
  U_{k-j-1}(q)(aU_{i-n}(q)+bT_{i-n}(q))
  +U_{j-i-1}(q)&(aU_{k-n}(q)+bT_{k-n}(q))\\
  &=U_{k-i-1}(q)(aU_{j-n}(q)+bT_{j-n}(q)),
}
which, in view of \eq{pi}, shows that \eq{ijk=} holds.

Finally, let $i\in\{n,\dots,m\}$ and $j\in\{1,\dots,\min(i-n,m-i)\}$. In view of \eq{UT}, we have that 
\Eq{*}{
  U_{i-j-n}(q)+U_{i+j-n}(q)=2T_j(q)U_{i-n}(q), \qquad
  T_{i-j-n}(q)+T_{i+j-n}(q)=2T_j(q)T_{i-n}(q).
}
Multiplying the first and second equalities by $a$ and $b$, respectively, and then adding them up side by side, we obtain
\Eq{*}{
  p_{i-j}+p_{i+j}
  &=(aU_{i-j-n}(q)+bT_{i-j-n}(q))
  +(aU_{i+j-n}(q)+bT_{i+j-n}(q))\\
  &=2T_j(q)(aU_{i-n}(q)+bT_{i-n}(q))
  =2T_j(q)p_i.
}
This completes the proof of \eq{ijT}.
\end{proof}

In the following statement, we establish some properties of the class of $q$-concave (and hence of $q$-convex) sequences.

\Prp{Props}{The cone $\C_q^\cap(n|m)$ is closed with respect to the pointwise minimum and the cone $\C_q^\cup(n|m)$ is closed with respect to the pointwise maximum.}

\begin{proof} To prove the statement for $\C_q^\cap(n|m)$, let $p,r\in\C_q^\cap(n|m)$ be arbitrary and denote $s:=\min(p,r)$ (i.e., $s_i:=\min(p_i,r_i)$ for all $i\in\{n,\dots,m\}$). Let $i\in\{n+1,\dots,m-1\}$. Then, by the $q$-concavity of $p$ and $r$, we have
\Eq{*}{
  s_{i-1}+s_{i+1}\leq p_{i-1}+p_{i+1}\leq qp_i \qquad\mbox{and}\qquad
  s_{i-1}+s_{i+1}\leq r_{i-1}+r_{i+1}\leq qr_i.
}
Therefore,
\Eq{*}{
  s_{i-1}+s_{i+1}\leq\min(qp_i,qr_i)=q\min(p_i,r_i)=qs_i,
}
which shows that $s$ is also $q$-concave. The proof of the statement for $\C_q^\cup(n|m)$ is analogous.
\end{proof}

As $q$-affine sequences are $q$-concave and also $q$-convex, we obtain that the pointwise minimum and maximum of a finite family of $q$-affine sequences are $q$-concave and also $q$-convex, respectively.

\Prp{qJC}{Let $i,j,k\in\{n,\dots,m\}$ with $i<j<k$. Assume that
\Eq{qnm}{
  q\geq \cos\bigg(\frac{\pi}{\max(j-i,k-j)}\bigg).
}
Then, for all $p\in\C_q^\cap(n|m)$,
\Eq{ijk}{
  U_{k-j-1}(q)p_{i}+U_{j-i-1}(q)p_{k}\leq U_{k-i-1}(q)p_j.
}
In particular, if $i\in\{n+1,\dots,m-1\}$ and $j\in\{1,\dots,\min(i-n,m-i)\}$ and
\Eq{qj0}{
  q>\cos\bigg(\frac{\pi}{j}\bigg),
}
then
\Eq{ijT+}{
  p_{i-j}+p_{i+j}\leq 2T_{j}(q)p_i.
} 
}

\begin{proof} We shall verify \eq{ijk} by induction
on $\ell:=k-i$. If $\ell=2$, that is, $j-i=k-j=1$, then \eq{ijk} is equivalent to the $q$-concavity of $p$, because $U_0(q)=1$ and $U_1(q)=2q$.

Assume that we have verified \eq{ijk} for all $i<j<k$ with $k-i\leq\ell$, where $\ell\geq2$. Suppose that $k-i=\ell+1\geq3$ and \eq{qnm} holds. Then $\max(j-i,k-j)\geq2$. We now distinguish two cases.

The first the case is when $j-i\geq2$. Then $k-(i+1)=\ell$ and $j-i\leq k-i-1\leq\ell$ and, using \eq{qnm}, it follows that
\Eq{*}{
  q\geq \cos\bigg(\frac{\pi}{\max(j-(i+1),k-j)}\bigg)
  \qquad\mbox{and}\qquad 
  q\geq \cos\bigg(\frac{\pi}{\max((i+1)-i,j-(i+1))}\bigg).
}
Thus, applying the inductive hypotheses for the triplets $i+1<j<k$ and for $i<i+1<j$, we obtain
\Eq{*}{
  U_{k-j-1}(q)p_{i+1}+U_{j-i-2}(q)p_k&\leq U_{k-i-2}(q)p_j,\\
  U_{j-i-2}(q)p_{i}+U_{0}(q)p_{j}&\leq U_{j-i-1}(q)p_{i+1}.
}
The inequality \eq{qnm} shows that $q$ is nonsmaller than the largest roots of $U_{j-i-1}$ and $U_{k-j-1}$, hence $U_{j-i-1}(q)\geq0$ and $U_{k-j-1}(q)\geq0$. Multiplying the first inequality by $U_{j-i-1}(q)$, the second one by $U_{k-j-1}(q)$, and adding up the inequalities so obtained side by side, we get
\Eq{*}{
  U_{k-j-1}(q)U_{j-i-2}(q)p_{i}+U_{j-i-1}(q)U_{j-i-2}(q)p_{k}
  \leq \big(U_{j-i-1}(q)U_{k-i-2}(q)-U_{k-j-1}(q)U_{0}(q)\big)p_j.
}
On the other hand, applying \lem{IdU} for the numbers $k-j-1<k-i-2<k-i-1$, we have that
\Eq{*}{
 U_{j-i-1}(q)U_{k-i-2}(q)
 =U_{0}(q)U_{k-j-1}(q)+U_{j-i-2}(q)U_{k-i-1}(q). 
}
Therefore, the above inequality can be rewritten as
\Eq{*}{
  U_{k-j-1}(q)U_{j-i-2}(q)p_{i}+U_{j-i-1}(q)U_{j-i-2}(q)p_{k}
  \leq U_{j-i-2}(q)U_{k-i-1}(q)p_j.
}
By \eq{qnm}, $q$ is strictly bigger than $\cos\big(\frac{\pi}{j-i-1}\big)$, which is the largest root of $U_{j-i-2}$ if $j-i>2$, therefore $U_{j-i-2}(q)>0$. If $i-j=2$, then $U_{j-i-2}(q)=U_0(q)=1>0$. Now dividing the last inequality by this positive value side by side, we arrive at the desired inequality \eq{ijk}. 

The proof in the second case when $k-j\geq2$ is completely analogous, therefore it is omitted. 

Finally, let $i\in\{n+1,\dots,m-1\}$ and $j\in\{1,\dots,\min(i-n,m-i)\}$ and assume that \eq{qj0} is satisfied. We apply the previous statement to the triplet $(i-j,i,i+j)$. Then, also using identity \eq{UT}, we get
\Eq{ijU}{
  U_{j-1}(q)p_{i-j}+U_{j-1}(q)p_{i+j}
  \leq U_{2j-1}(q)p_i=2U_{j-1}(q)T_j(q)p_i.
} 
In view of \eq{qj0}, we have that $q$ is bigger than the largest root of $U_{j-1}$ if $j\geq 2$, hence $U_{j-1}(q)>0$. This inequality is obviously true if $j=1$. Thus, after dividing \eq{ijU} by $U_{j-1}(q)$ side by side, this inequality implies \eq{ijT+}.
\end{proof}

\Prp{3}{Let $j,k\in\{n,\dots,m\}$ with $j<k$. In addition, assume that
\Eq{qnm+}{
  q>\cos\Big(\frac{\pi}{k-j}\Big).
}
Let $p\in \C^\cap(n|m)$ and define
\Eq{*}{
  r_i:=p_k\frac{U_{i-j-1}(q)}{U_{k-j-1}(q)}+p_j\frac{U_{k-i-1}(q)}{U_{k-j-1}(q)} \qquad(i\in \{n,\dots,m\}).
}
Then, $r=(r_n,\dots,r_m)$ is a $q$-affine sequence and, for $i\in\{n,\dots,m\}$,
\Eq{*}{
  r_i\begin{cases}
     \geq p_i & \text{ if }\  i<j \text{ or } k<i.\\
     =p_i & \text{ if }\ i\in\{j,k\}. \\
     \leq p_i & \text{ if }\  j<i<k.
     \end{cases}
}
}

\begin{proof} If $k-j=1$, then $U_{k-j-1}(q)=U_0(q)=1>0$. If $k-j\geq2$, then $q$ is bigger than the largest root of $U_{k-j-1}$. Therefore $U_{k-j-1}(q)>0$ and hence the sequence $(r_i)$ is well-defined. From the recursive formula \eq{T} of Chebyshev polynomials of the second kind, for $i\in\{n+1,\dots,m-1\}$, it follows that
\Eq{*}{
  U_{(i-1)-j-1}(q)+U_{(i+1)-j-1}(q)=2qU_{i-j-1}(q), \qquad
  U_{k-(i-1)-1}(q)+U_{k-(i+1)-1}(q)=2qU_{k-i-1}(q).
}
Multiplying theses equalities by $\frac{p_k}{U_{k-j-1}(q)}$ and by 
$\frac{p_j}{U_{k-j-1}(q)}$, respectively, and then adding them up side by side, we obtain that $r_{i-1}+r_{i+1}=2qr_i$, which shows that $(r_i)$ is a $q$-affine sequence.

If $i=j$, or $i=k$, then, by $U_{-1}=0$, we can see that $r_j=p_j$ and $r_k=p_k$. Suppose first that $j<i<k$. From the equality \eq{ijk=} of the second assertion of \prp{qJA} applied to the $q$-affine sequence $(r_i)$, we get
\Eq{*}{
  U_{k-j-1}(q)r_i=U_{k-i-1}(q)r_j+U_{i-j-1}(q)r_k.
}
On the other hand, applying inequality \eq{ijk} of \prp{qJC} for the $q$-concave sequence $(p_i)$, we get
\Eq{*}{
 U_{k-i-1}(q)p_j+U_{i-j-1}(q)p_k\leq U_{k-j-1}(q)p_i
}
Using that $r_j=p_j$ and $r_k=p_k$, it follows that
\Eq{*}{
  U_{k-j-1}(q)r_i
  =U_{k-i-1}(q)r_j+U_{i-j-1}(q)r_k
  =U_{k-i-1}(q)p_j+U_{i-j-1}(q)p_k\leq U_{k-j-1}(q)p_i,
}
which, by $U_{k-j-1}(q)>0$ simplifies to the inequality $r_i\leq p_i$.

For the remaining inequalities, suppose first that $i<j$. By the $q$ affinity of $(r_i)$, the second assertion of \prp{qJA} implies
\Eq{*}{
  U_{k-i-1}(q)r_j=U_{k-j-1}(q)r_i+U_{j-i-1}(q)r_k
}
and hence
\Eq{*}{
  U_{k-j-1}(q)r_i=U_{k-i-1}(q)r_j-U_{j-i-1}(q)r_k.
}
On the other hand, applying inequality \eq{ijk} of \prp{qJC} for the $q$-concave sequence $(p_i)$, we get
\Eq{*}{
 U_{k-j-1}(q)p_i+U_{j-i-1}(q)p_k\leq U_{k-i-1}(q)p_j
}
and hence
\Eq{*}{
 U_{k-i-1}(q)p_j-U_{j-i-1}(q)p_k\geq U_{k-j-1}(q)p_i.
}
Combining these inequalities and using $r_j=p_j$ and $r_k=p_k$, we can conclude that
\Eq{*}{
  U_{k-j-1}(q)r_i=U_{k-i-1}(q)r_j-U_{j-i-1}(q)r_k
  =U_{k-i-1}(q)p_j-U_{j-i-1}(q)p_k\geq U_{k-j-1}(q)p_i.
}
This inequality, by $U_{k-j-1}(q)>0$, is equivalent to $r_i\geq p_i$ as desired.

The proof of $r_i\geq p_i$ in the case $k<i$ is completely similar and therefore omitted.
\end{proof}

In the following proposition, we establish a characterization of $q$-concave sequences.

\Prp{4}{Let $p\in\S(n|m)$. Then $p$ is $q$-concave if and only if, for all $j\in\{n,\dots,m-1\}$, there exists $r\in\A_q(n|m)$ such that 
\Eq{pr}{
p_j=r_j, \qquad p_{j+1}=r_{j+1}, \qquad\mbox{and}\qquad p_i\leq r_i \quad\mbox{for}\quad i\in\{n,\dots,m\}.
}}

\begin{proof} Assume first that $p$ is $q$-concave and let $j\in\{n,\dots,m-1\}$. Then, with $k=j+1$, we can see that \eq{qnm+} holds, therefore applying \prp{3}, the sequence $r\in\S(n|m)$ defined by
\Eq{*}{
  r_i:=p_{j+1}U_{i-j-1}(q)+p_jU_{j-i}(q)
}
is $q$-affine and satisfies all te conditions in \eq{pr}.

To prove the sufficiency part of the assertion, assume that
$j\in\{n,\dots,m-1\}$, there exists $q$-affine sequence $r^j\in\A_q(n|m)$ such that 
\Eq{*}{
p_j=r_j^j, \qquad p_{j+1}=r_{j+1}^j, \qquad\mbox{and}\qquad p_i\leq r_i^j \quad\mbox{for}\quad i\in\{n,\dots,m\}.
}
Then, it follows that
\Eq{*}{
   p_i=\min_{n\leq j\leq m-1} r_i^j,
}
which shows that $p$ is the pointwise minimum of finitely many (in fact, $m-n$) $q$-affine sequences. Thus, by \prp{Props}, it follows that $p$ is $q$-concave.
\end{proof}

\section{A minimax-type problem}

Throughout this section, $n,m$ are integers with $2\leq m-n$ and we consider the following minimum problem: Let $M:\R_+^{m-n-1}\to\R_+$ be an $(m-n-1)$-variable mean. Our aim is to find the largest nonnegative constant $C_M$ such that, for all $p\in\S(n|m)$ with $p_n,p_m\geq0$ and $p_{n+1},\dots,p_{m-1}>0$,
\Eq{*}{
  C_M\leq M\bigg(\frac{p_{n}+p_{n+2}}{2p_{n+1}},\dots,
  \frac{p_{i-1}+p_{i+1}}{2p_i},\dots,\frac{p_{m-2}+p_{m}}{2p_{m-1}}\bigg).
}
By taking $p$ as a constant sequence, one can see that the right hand side of this inequality then equals $1$, hence it follows that $C_M\leq1$. As we shall see below, this estimate can be essentially improved for several concrete means.

In the case when $M$ is the $(m-n-1)$-variable arithmetic mean $A_{m-n-1}$, we can obtain the following result.

\Prp{AM}{$C_A=\frac{m-n-2}{m-n-1}$, that is, for all $p\in\S(n|m)$ with $p_n,p_m\geq0$ and $p_{n+1},\dots,p_{m-1}>0$,
\Eq{nm}{
  \frac{m-n-2}{m-n-1}
  \leq \frac{1}{m-n-1}\sum_{i=n+1}^{m-1}\frac{p_{i-1}+p_{i+1}}{2p_i}
}
and the constant on the left hand side is the best possible.}

\begin{proof} If $m-n=2$, that is, $m=n+2$, then the left hand side of \eq{nm} equals zero, thus, the inequality is trivial. On the other hand, for $(p_n,p_{n+1},p_{n+2})=(0,1,0)$ equality holds in \eq{nm}.
Thus, in the rest of the proof, we may assume that $m-n>2$. 

To prove \eq{nm}, let $p\in\S(n|m)$ with $p_n,p_m\geq0$ and $p_{n+1},\dots,p_{m-1}>0$. Then (using the arithmetic-geometric mean inequality in the last step), we obtain
\Eq{*}{
  \sum_{i=n+1}^{m-1}\frac{p_{i-1}+p_{i+1}}{2p_i}
  &=\frac{p_{n}+p_{n+2}}{2p_{n+1}}+\sum_{i=n+2}^{m-2}\frac{p_{i-1}+p_{i+1}}{2p_i}+\frac{p_{m-2}+p_{m}}{2p_{m-1}}\\
  &\geq\frac{p_{n+2}}{2p_{n+1}}+\sum_{i=n+2}^{m-2}\bigg(\frac{p_{i-1}}{2p_i}+\frac{p_{i+1}}{2p_i}\bigg)+\frac{p_{m-2}}{2p_{m-1}}\\
  &=\sum_{i=n+1}^{m-2}\frac12\bigg(\frac{p_{i+1}}{p_i}+\frac{p_i}{p_{i+1}}\bigg)
  \geq \sum_{i=n+1}^{m-2}\sqrt{\frac{p_{i+1}}{p_i}\cdot\frac{p_i}{p_{i+1}}}=m-n-2.
}
Dividing the above obtained inequality by $m-n-1$ side by side, we can see that \eq{nm} holds. On the other hand, for $(p_n,p_{n+1},\dots,p_{m-1},p_m)=(0,1,\dots,1,0)$ equality holds in \eq{nm}, therefore, the left hand side of \eq{nm} is the largest possible, indeed. 
\end{proof}

In order to reach a higher level of generality, for $r\in[-\infty,\infty]$ and $k\in\N$, we define the \emph{$k$-variable $r$th power mean (or Hölder mean)} of the variables $u_1,\dots,u_k\in\R_+$ by 
\Eq{*}{
  H_{r,k}(u_1,\dots,u_k):=
  \begin{cases} 
  \min(x_1,\dots,x_k) & \mbox{if } r=-\infty,\\[2mm]
  \bigg(\dfrac{u_1^r+\dots+u_k^r}{k}\bigg)^{\frac{1}{r}} 
  & \mbox{if } r\in\R\setminus\{0\},\\[3mm]
  \sqrt[k]{u_1\cdots u_k} & \mbox{if } r=0,\\[2mm]
  \max(x_1,\dots,x_k) & \mbox{if } r=\infty.
  \end{cases}
}
Obviously, the mean $H_{1,k}$ equals the $k$-variable arithmetic mean $A_k$ and $H_{0,k}$ equals the $k$-variable geometric mean $G_k$. It is well known that, for all $k\in\N$ and $-\infty\leq r\leq s\leq\infty$, the comparison inequality $H_{r,k}\leq H_{s,k}$ holds. In particular, $G_k\leq A_k$, which is the celebrated inequality between the geometric and arithmetic means. 

For the investigation of the more general problem in terms of power means, for $r\in\R$ and $k\in\N$, we introduce the function $F_{r,k}:\R_+^k\to\R_+$  by
\Eq{*}{
  F_{r,k}(u_1,\dots,u_k)
  :=u_1^r+\sum_{i=1}^{k-1}\Big(\frac{1}{u_{i}}+u_{i+1}\Big)^r+\frac1{u_k^r}.
}

\Lem{F}{Let $r>0$ and $k\in\N$. Then
\Eq{Frk}{
F_{r,k}\geq
\begin{cases}
2 &  \mbox{if } k=1,\\
2^{\frac{r+1}{2}}+(k-2)2^r+2^{\frac{r+1}{2}}
  & \mbox{if } k\geq 2,\, r\leq 1,\\
2^rk^{1-r}\big(2^{\frac{1-r}{2r}}+(k-2)+2^{\frac{1-r}{2r}}\big)^r 
  & \mbox{if } k\geq 2,\, r\geq 1,\\
\end{cases}
}
and the estimates are sharp if $k\in\{1,2\}$ or $r=1$. Furthermore, for all $k\in\N$
\Eq{Frk+}{
F_{r,k}\geq k2^{r+\frac{1-r}{k}},
}
which is also sharp if if $k\in\{1,2\}$ or $r=1$.
In the particular case when $k$ is odd, we also have that
\Eq{Frk++}{
  F_{r,k}\geq k+1,
}
which is sharp if $k=1$ and which is sharper than \eq{Frk} and \eq{Frk+} if $r$ is a sufficiently small positive number.}

\begin{proof} 
If $k=1$ then, by the arithmetic-geometric mean inequality, for all $u_1\in\R_+$, we easily get  
\Eq{*}{
  F_{r,1}(u_1)=u_1^r+\frac1{u_1^r}
  =2A_2\Big(u_1^r,\frac1{u_1^r}\Big)
  \geq 2G_2\Big(u_1^r,\frac1{u_1^r}\Big)
  =2\sqrt{u_1^r\cdot\frac1{u_1^r}}=2.
}
Observe that $F_{r,1}(1)=2$, hence the lower estimate $2$ is best possible in this case.

Now assume that $r\leq 1$ and $k\geq2$, and let $u_1,\dots,u_k\in\R_+$ be arbitrary. Then, by the comparison inequality $H_{r,2}\leq H_{1,2}=A_2$, for all $i\in\{1,\dots,k-1\}$, we get
\Eq{*}{
  \frac12\Big(\frac{1}{u_{i}}+u_{i+1}\Big)
  =A_2\Big(\frac{1}{u_{i}},u_{i+1}\Big)
  \geq H_{r,2}\Big(\frac{1}{u_{i}},u_{i+1}\Big)
  =\bigg(\frac12\Big(\frac{1}{u_i^r}+(u_{i+1})^r\Big)\bigg)^{\frac1r}.
}
Using this inequality and arithmetic-geometric mean inequality at the end, we obtain
\Eq{*}{
  F_{r,k}(u_1,\dots,u_k)
  &:=u_1^r+\sum_{i=1}^{k-1}2^r\Big(\frac12\Big(\frac{1}{u_{i}}+u_{i+1}\Big)\Big)^r+\frac1{u_k^r}\\
  &\geq u_1^r+\sum_{i=1}^{k-1}2^r\frac12\Big(\frac{1}{u_i^r}+(u_{i+1})^r\Big)+\frac1{u_k^r}\\
  &=u_1^r+2^{r-1}\frac{1}{u_1^r}+\sum_{i=2}^{k-1}2^{r-1}\Big(\frac{1}{u_i^r}+u_i^r\Big)+2^{r-1}u_k+\frac1{u_k}\\
  &=2A_2\Big(u_1^r,2^{r-1}\frac{1}{u_1^r}\Big)
  +\sum_{i=2}^{k-1}2^{r}A_2\Big(\frac{1}{u_i^r},u_i^r\Big)
  +2A_2\big(2^{r-1}u_k,\frac1{u_k}\Big)\\
  &\geq2G_2\Big(u_1^r,2^{r-1}\frac{1}{u_1^r}\Big)
  +\sum_{i=2}^{k-1}2^{r}G_2\Big(\frac{1}{u_i^r},u_i^r\Big)
  +2G_2\big(2^{r-1}u_k,\frac1{u_k}\Big)\\
  &=2^{\frac{r+1}{2}}+(k-2)2^r+2^{\frac{r+1}{2}}.
}
This proves the assertion when $r\leq 1$ and $k\geq2$. 
Finally, by arithmetic-geometric mean inequality again, we get
\Eq{*}{
  F_{r,k}(u_1,\dots,u_k)
  &\geq2^{\frac{r+1}{2}}+(k-2)2^r+2^{\frac{r+1}{2}}
  =kA_k\big(2^{\frac{r+1}{2}},2^r,\dots,2^r,2^{\frac{r+1}{2}}\big)\\
  &\geq kG_k\big(2^{\frac{r+1}{2}},2^r,\dots,2^r,2^{\frac{r+1}{2}}\big)
  =k\big(2^{(k-2)r+r+1}\big)^{\frac1k}=k2^{r+\frac{1-r}{k}},
}
which shows that \eq{Frk+} is also valid.

In the case $r\geq 1$ and $k\geq2$, using the comparison inequality $A_{2k}=H_{1,2k}\geq H_{\frac1r,2k}$ and the $2$-variable arithmetic-geometric mean inequality, we obtain
\Eq{*}{
  F_{r,k}(u_1,\dots,u_k)
  &=u_1^r+\sum_{i=1}^{k-1}2\cdot\frac12\Big(\frac{1}{u_{i}}+u_{i+1}\Big)^r+\frac1{u_k^r}\\
  &=2kA_{2k}\bigg(u_1^r,\dots,\frac12\Big(\frac{1}{u_{i}}+u_{i+1}\Big)^r,\frac12\Big(\frac{1}{u_{i}}+u_{i+1}\Big)^r,\dots\frac1{u_k^r}\bigg)\\
  &\geq 2kH_{\frac1r,2k}\bigg(u_1^r,\dots,\frac12\Big(\frac{1}{u_{i}}+u_{i+1}\Big)^r,\frac12\Big(\frac{1}{u_{i}}+u_{i+1}\Big)^r,\dots\frac1{u_k^r}\bigg)\\
  &=2k\bigg(\frac1{2k}\bigg(u_1+\sum_{i=1}^{k-1}2\cdot 2^{-\frac1r}\Big(\frac{1}{u_{i}}+u_{i+1}\Big)+\frac1{u_k}\bigg)\bigg)^r\\
  &= (2k)^{1-r}\bigg(u_1+2^{1-\frac1r}\frac{1}{u_{1}}+\sum_{i=2}^{k-1}2^{1-\frac1r}\Big(\frac{1}{u_{i}}+u_{i}\Big)+2^{1-\frac1r}u_k+\frac1{u_k}\bigg)^r\\
  &\geq (2k)^{1-r}\big(2\cdot 2^{\frac{r-1}{2r}}+2(k-2)2^{1-\frac1r}+2\cdot 2^{\frac{r-1}{2r}}\big)^r\\
  &=2^rk^{1-r}\big(2^{\frac{1-r}{2r}}+(k-2)+2^{\frac{1-r}{2r}}\big)^r.
}
This proves the assertion when $r\geq 1$ and $k\geq2$. 
Finally, by arithmetic-geometric mean inequality again, we get
\Eq{*}{
  F_{r,k}(u_1,\dots,u_k)
  &\geq2^rk^{1-r}\big(2^{\frac{1-r}{2r}}+(k-2)+2^{\frac{1-r}{2r}}\big)^r
  =2^rkA_k\big(2^{\frac{1-r}{2r}},1,\dots,1,2^{\frac{1-r}{2r}}\big)^r\\
  &\geq 2^rkG_k\big(2^{\frac{1-r}{2r}},1,\dots,1,2^{\frac{1-r}{2r}}\big)^r
  =2^rk \sqrt[k]{2^{\frac{1-r}{r}}}=k2^{r+\frac{1-r}{k}},
}
which shows that \eq{Frk+} is also valid.

If $k=2$, then the lower estimates \eq{Frk} and \eq{Frk+} simplify to the inequality
\Eq{*}{
 F_{r,2}\geq 2^{\frac{3+r}{2}}.
}
On the other hand, with $u_1:=2^{\frac{r-1}{2r}}$ and $u_2:=2^{\frac{1-r}{2r}}$, one can see that
\Eq{*}{
 F_{r,2}(u_1,u_2)=u_1^r+\Big(\frac{1}{u_{1}}+u_{2}\Big)^r+\frac1{u_2^r}
  =2^{\frac{r+1}{2}}+2^{\frac{1+r}{2}}=2^{\frac{3+r}{2}},
}
which proves that the lower estimate $2^{\frac{3+r}{2}}$ is sharp.

If $r=1$, then all the lower estimates simplify to the inequality
\Eq{*}{
 F_{1,k}\geq 2k,
}
which is attained at $u_1=\dots=u_k=1$. This proves that the lower estimate $2k$ is sharp in this case.

Finally, we prove that \eq{Frk++} holds. This inequality is a consequence of \eq{Frk} in the case $k=1$. Thus, we may assume that $k\geq3$ is odd. Then, for $u_1,\dots,u_k\in\R_+$, we get
\Eq{*}{
  F_{r,k}(u_1,\dots,u_k)
  &=u_1^r+\Big(\frac{1}{u_{1}}+u_{2}\Big)^r
  +\sum_{i=2}^{k-2}\Big(\frac{1}{u_{i}}+u_{i+1}\Big)^r
  +\Big(\frac{1}{u_{k-1}}+u_{k}\Big)^r+\frac1{u_k^r}\\
  &\geq u_1^r+\frac{1}{u_{1}^r}
  +\sum_{j=1}^{\frac{k-3}2}\bigg(\Big(\frac{1}{u_{2j}}+u_{2j+1}\Big)^r
  +\Big(\frac{1}{u_{2j+1}}+u_{2j+2}\Big)^r\bigg)
  +u_{k}^r+\frac1{u_k^r}\\
  &\geq u_1^r+\frac{1}{u_{1}^r}
  +\sum_{j=1}^{\frac{k-3}2}\bigg(u_{2j+1}^r+\frac{1}{u_{2j+1}^r}\bigg)
  +u_{k}^r+\frac1{u_k^r}\\
  &=\sum_{j=0}^{\frac{k-1}2}\bigg(u_{2j+1}^r+\frac{1}{u_{2j+1}^r}\bigg)
  \geq \sum_{j=0}^{\frac{k-1}2}2=k+1.
}
If $r$ tends to zero in \eq{Frk}, then the limit of the lower estimate is $2\sqrt{2}+k-2$, which is smaller than $k+1$, showing that \eq{Frk++} provides a better lower estimate than \eq{Frk} for small positive values of $r$.
\end{proof}

\Prp{PM}{Let $r>0$. Then
\Eq{Hr}{
  C_{H_{r,m-n-1}}\geq
  \begin{cases}
   \dfrac{1}{2} &\mbox{if } m=n+3, \\[3mm]
   \bigg(\dfrac{2^{\frac{1-r}{2}}+(m-n-4)+2^{\frac{1-r}{2}}}{m-n-1}\bigg)^{\frac{1}{r}} &\mbox{if } m\geq n+4,\, 0<r\leq1, \\[4mm]
   \Big(\dfrac{m-n-2}{m-n-1}\Big)^{\frac{1}{r}}\cdot
   \dfrac{2^{\frac{1-r}{2r}}+(m-n-4)+2^{\frac{1-r}{2r}}}{m-n-2}&\mbox{if } m\geq n+4,\, 1\leq r. \\
  \end{cases}
}
and the constant on the left hand side is the best possible if either $m\in\{n+3,n+4\}$ or $r=1$. In addition, if $m-n$ is odd, then
\Eq{Hr+}{
  C_{H_{r,m-n-1}}\geq\frac{1}{2}.
}}

\begin{proof} If $m-n=2$, that is, $m=n+2$, then the left hand side of \eq{nm} equals zero, thus, the inequality is trivial. On the other hand, for $(p_n,p_{n+1},p_{n+2})=(0,1,0)$ equality holds in \eq{nm}.
Thus, in the rest of the proof, we may assume that $m-n>2$. 

To prove \eq{nm}, let $p\in\S(n|m)$ with $p_n,p_m\geq0$ and $p_{n+1},\dots,p_{m-1}>0$. Then
\Eq{*}{
  2^r\sum_{i=n+1}^{m-1}\Big(\frac{p_{i-1}+p_{i+1}}{2p_i}\Big)^r
  &=\Big(\frac{p_{n}+p_{n+2}}{p_{n+1}}\Big)^r+\sum_{i=n+2}^{m-2}\Big(\frac{p_{i-1}+p_{i+1}}{p_i}\Big)^r+\Big(\frac{p_{m-2}+p_{m}}{p_{m-1}}\Big)^r\\
  &\geq\Big(\frac{p_{n+2}}{p_{n+1}}\Big)^r+\sum_{i=n+2}^{m-2}\Big(\frac{p_{i-1}}{p_i}+\frac{p_{i+1}}{p_i}\Big)^r+\Big(\frac{p_{m-2}}{p_{m-1}}\Big)^r\\
  &=F_{r,m-n-2}\bigg(\frac{p_{n+2}}{p_{n+1}},\dots,\frac{p_{m-3}}{p_{m-2}}\bigg).
}
Therefore,
\Eq{*}{
  \bigg(\frac{1}{m-n-1}\sum_{i=n+1}^{m-1}\Big(\frac{p_{i-1}+p_{i+1}}{2p_i}\Big)^r\bigg)^{\frac{1}{r}}
  &\geq \bigg(\frac{1}{2^r(m-n-1)}F_{r,m-n-2}\bigg(\frac{p_{n+2}}{p_{n+1}},\dots,\frac{p_{m-3}}{p_{m-2}}\bigg)\bigg)^{\frac{1}{r}}.
}
If $m=n+3$, then, by the $k=1$ case of \lem{F}, we get
\Eq{*}{
  \bigg(\frac{1}{2}\sum_{i=n+1}^{n+2}\Big(\frac{p_{i-1}+p_{i+1}}{2p_i}\Big)^r\bigg)^{\frac{1}{r}}\geq \frac{1}{2}.
}
Applying \lem{F}, for $k:=m-n-2\geq 2$ and $0<r\leq 1$, we get
\Eq{*}{
  \bigg(\frac{1}{m-n-1}\sum_{i=n+1}^{m-1}\Big(\frac{p_{i-1}+p_{i+1}}{2p_i}\Big)^r\bigg)^{\frac{1}{r}}
  \geq\bigg(\frac{2^{\frac{1-r}{2}}+(m-n-4)+2^{\frac{1-r}{2}}}{m-n-1}\bigg)^{\frac{1}{r}}.
}
Similarly, for $k:=m-n-2\geq 2$ and $r\geq 1$, it follows that 
\Eq{*}{
  \bigg(\frac{1}{m-n-1}\sum_{i=n+1}^{m-1}\Big(\frac{p_{i-1}+p_{i+1}}{2p_i}\Big)^r\bigg)^{\frac{1}{r}}
  \geq\Big(\frac{m-n-2}{m-n-1}\Big)^{\frac{1}{r}}\cdot
  \frac{2^{\frac{1-r}{2r}}+(m-n-4)+2^{\frac{1-r}{2r}}}{m-n-2}.
}
To prove \eq{Hr+}, assume that $m-n$ is odd. Then, applying the inequality \eq{Frk++} for $k=m-n-2$, we get
\Eq{*}{
  \bigg(\frac{1}{m-n-1}\sum_{i=n+1}^{m-1}\Big(\frac{p_{i-1}+p_{i+1}}{2p_i}\Big)^r\bigg)^{\frac{1}{r}}
  \geq\bigg(\frac{(m-n-2)+1}{2^r(m-n-1)}\bigg)^{\frac{1}{r}}=\frac12,
}
which was to be shown.
\end{proof}

In the case when $M$ is the $(m-n-1)$-variable geometric mean $G$, we can establish the following result in which we will get an exact formula for the constant $C_{G_k}$.

\Prp{GM}{$C_{G_{m-n-1}}=\frac{1+(-1)^{m-n-1}}{4}$, that is, for all sequences $p\in\S(n|m)$ with $p_n,p_m\geq0$ and $p_{n+1},\dots,p_{m-1}>0$,
\Eq{nmG}{
  \frac{1+(-1)^{m-n-1}}{4}\leq\sqrt[m-n-1]{\prod_{i=n+1}^{m-1}\frac{p_{i-1}+p_{i+1}}{2p_i}}.
}
and the constant on the left hand side is the best possible.}

\begin{proof} Assume first that $m-n$ is even. Then the left hand side of \eq{nmG} equals zero, thus, the inequality is trivial. To show that the left hand side is optimal, define the sequence $p\in\S(n|m)$ by
\Eq{*}{
p_{n+2i}:=\varepsilon \qquad(i\in\{0,\dots,\tfrac{m-n}{2}\})
\qquad\mbox{and}\qquad p_{n+2i+1}:=1 \qquad(i\in\{0,\dots,\tfrac{m-n-2}{2}\})
}
where $\varepsilon>0$ is an arbitrary positive number.
Then, using that $m-n$ is even, we can obtain that
\Eq{*}{
  \prod_{i=n+1}^{m-1}\frac{p_{i-1}+p_{i+1}}{2p_i}
  =\prod_{i=n+1}^{m-1}\varepsilon^{(-1)^{i-(n+1)}}=\varepsilon.
}
Therefore, the rights hand side of \eq{nmG} equals $\sqrt[m-n-1]{\varepsilon}$, which can be arbitrarily small. Hence, in this case, we obtain that $C_G=0$.

Consider now the case when $m-n$ is odd and $m-n\geq 3$. Using that the product has an even number of factors, we get
\Eq{*}{
  \prod_{i=n+1}^{m-1}\frac{p_{i-1}+p_{i+1}}{2p_i}
  &=\prod_{j=0}^{\frac{m-n-3}{2}}\frac{p_{n+2j}+p_{n+2+2j}}{2p_{n+1+2j}}
  \cdot\frac{p_{n+1+2j}+p_{n+3+2j}}{2p_{n+2+2j}}\\
  &\geq\prod_{j=0}^{\frac{m-n-3}{2}}\frac{p_{n+2+2j}}{2p_{n+1+2j}}
  \cdot\frac{p_{n+1+2j}}{2p_{n+2+2j}}=\frac1{2^{m-n-1}}.
}
Taking the $(m-n-1)$th root of this inequality side by side, we obtain that \eq{nmG} is also true in the case when $m-n$ is odd and $m-n\geq 3$.

To verify the sharpness of the left hand side of \eq{nmG}, let $\varepsilon>0$ be  arbitrary and, for $i\in\{n,\dots,m\}$, define
\Eq{*}{
  p_i:=\begin{cases}
       \varepsilon^{\frac{m-i-1}{2}}&\mbox{if } i-n \mbox{ is even,}\\[2mm]
       \varepsilon^{\frac{i-n-1}{2}}&\mbox{if } i-n \mbox{ is odd.}
       \end{cases}
}
Then
\Eq{*}{
  \prod_{i=n+1}^{m-1}\frac{p_{i-1}+p_{i+1}}{2p_i}
  &=\prod_{j=0}^{\frac{m-n-3}{2}}\frac{p_{n+2j}+p_{n+2+2j}}{2p_{n+1+2j}}\cdot\frac{p_{n+1+2j}+p_{n+3+2j}}{2p_{n+2+2j}}\\
  &=\prod_{j=0}^{\frac{m-n-3}{2}}\frac{\varepsilon^{\frac{m-n-2j-1}{2}}+\varepsilon^{\frac{m-n-2j-3}{2}}}{2\varepsilon^{j}}\cdot\frac{\varepsilon^{j}+\varepsilon^{j+1}}{2\varepsilon^{\frac{m-n-2j-3}{2}}}\\
  &=\prod_{j=0}^{\frac{m-n-3}{2}}\bigg(\frac{\varepsilon+1}{2}\cdot\frac{1+\varepsilon}{2}\bigg)=\bigg(\frac{1+\varepsilon}{2}\bigg)^{m-n-1}.
}
By taking $\varepsilon$ arbitrarily small, we can see that the right hand side of the above equality can be arbitrarily close to $\frac1{2^{m-n-1}}$, which shows that the left hand side of \eq{nmG} is a sharp lower bound for the right hand side.
\end{proof}

\Prp{Ch}{$C_{H_{\infty,m-n-1}}=\cos\big(\frac{\pi}{m-n}\big)$, that is, for all sequences $p\in\S(n|m)$ with $p_n,p_m\geq0$ and $p_{n+1},\dots,p_{m-1}>0$,  
\Eq{q1}{
  \cos\Big(\frac{\pi}{m-n}\Big)
  \leq\max_{n+1\leq i\leq m-1}\frac{p_{i-1}+p_{i+1}}{2p_i}.
}
Moreover, with $p_i:=\sin\big(\frac{i-n}{m-n}\pi\big)$, the inequality \eq{q1} holds with equality.
}

\begin{proof} Let 
\Eq{*}{
  q:=\max_{n+1\leq i\leq m-1}\frac{p_{i-1}+p_{i+1}}{2p_i}.
}
Then, using the positivity of $p_1,\dots,p_n$, it follows that the sequence $p$ is $q$-concave.

In the first part of the proof, we show that, for $k\in\{n,\dots,m-1\}$, 
\Eq{Sk}{
  0\leq U_{k-n}(q) \qquad\mbox{and}\qquad U_{k-n-1}(q)p_{k+1}\leq U_{k-n}(q)p_k.
}
These inequalities are obvious for $k=n$ because $U_0(q)=1$ and $U_{-1}(q)=0\leq p_n$. Assume that we have proved \eq{Sk} for some $k \in\{n,\dots,m-2\}$. Then, by the $q$-concavity of $p$, we have that
\Eq{*}{
  p_k+p_{k+2}\leq 2q p_{k+1} 
}
Multiplying this inequality by $U_{k-n}(q)\geq0$ and adding it to the second inequality in \eq{Sk} side by side, we get
\Eq{*}{
 U_{k-n-1}(q)p_{k+1}+U_{k-n}(q)p_{k+2}\leq 2q U_{k-n}(q)p_{k+1},
} 
which, by applying \eq{T}, implies
\Eq{*}{
  U_{k-n}(q)p_{k+2}\leq (2q U_{k-n}(q)-U_{k-n-1}(q))p_{k+1}=U_{k-n+1}(q)p_{k+1}.
}
This inequality shows that $U_{k-n+1}(q)$ is nonnegative and the second inequality in \eq{Sk} is valid for $k+1$ (instead of $k$).

Based on the first inequality in \eq{Sk}, for $k\in\{n,\dots,m-1\}$, we now show that
\Eq{Sk+}{
  \cos\Big(\frac{\pi}{k+1-n}\Big)\leq q.
}
This is obvious if $k=n$, since $q$ is nonnegative. If $k=n+1$, then \eq{Sk} gives that $0\leq U_1(q)=2q$ and hence $q\geq0=\cos\big(\frac{\pi}{2}\big)$, which proves \eq{Sk+} in this case. Now assume that \eq{Sk+} holds for some
$k\in\{n+1,\dots,m-2\}$. The two largest zeroes of $U_{k+1-n}$ are $\cos\big(\frac{2\pi}{k+2-n}\big)$ and $\cos\big(\frac{\pi}{k+2-n}\big)$, furthermore $U_{k+1-n}(t)<0$ if $\cos\big(\frac{2\pi}{k+2-n}\big)<t<\cos\big(\frac{\pi}{k+2-n}\big)$ and $U_{k+1-n}(t)\geq0$ if $t\geq\cos\big(\frac{\pi}{k+2-n}\big)$. Observe that
$\frac{\pi}{k+2-n}<\frac{\pi}{k+1-n}<\frac{2\pi}{k+2-n}$. Therefore, 
$\cos\big(\frac{2\pi}{k+2-n}\big)<\cos\big(\frac{\pi}{k+1-n}\big) <\cos\big(\frac{\pi}{k+2-n}\big)$. If $q$ were smaller than $\cos\big(\frac{\pi}{k+2-n}\big)$, then, by the inductive assumption, $\cos\big(\frac{\pi}{k+1-n}\big)\leq q <\cos\big(\frac{\pi}{k+2-n}\big)$ and hence $U_{k+1-n}(q)<0$, which contradicts \eq{Sk} (if it is applied for $k+1$ instead of $k$. Thus must be nonsmaller than $\cos\big(\frac{\pi}{k+2-n}\big)$, which shows that \eq{Sk+} is valid for $k+1$. 

Finally, applyinq \eq{Sk+} for $k=m-1$, we can conclude that $\cos\big(\frac{\pi}{m-n}\big)\leq q$, which proves that \eq{q1} holds.

To verify that \eq{q1} is sharp, let $p_i:=\sin\big(\frac{i-n}{m-n}\pi\big)$ for $i\in\{n,\dots,m\}$. Then, for $i\in\{n+1,\dots,m-1\}$, 
\Eq{*}{
  \frac{p_{i-1}+p_{i+1}}{2p_i}
  =\frac{\sin\big(\frac{i-1-n}{m-n}\pi\big)+\sin\big(\frac{i+1-n}{m-n}\pi\big)}{2\sin\big(\frac{i-n}{m-n}\pi\big)}
  =\frac{2\sin\big(\frac{i-n}{m-n}\pi\big)\cos\big(\frac{\pi}{m-n}\big)}{2\sin\big(\frac{i-n}{m-n}\pi\big)}
  =\cos\Big(\frac{\pi}{m-n}\Big),
}
which shows that \eq{q1} holds with equality for this particular sequence $p$.
\end{proof}

As a curiosity, we can obtain the following inequality for the cosine function.

\Cor{cos}{For $m\geq3$,
\Eq{cos}{
   \dfrac{m-4+\sqrt{2}}{m-2}\leq \cos\Big(\frac{\pi}{m}\Big),
}
and equality holds if $m=4$.}

\begin{proof} If $m=3$, then the inequality is equivalent to $\sqrt{2}-1\leq\frac12$, which is obviously true.

If $m\geq 4$ and $r\geq 1$, then, in view of \prp{PM}, for all sequences $p\in\S(0|m)$ with $p_0,p_m\geq0$ and $p_{1},\dots,p_{m-1}>0$, we have that
\Eq{*}{
   &\Big(\dfrac{m-2}{m-1}\Big)^{\frac{1}{r}}\cdot
   \dfrac{2^{\frac{1-r}{2r}}+(m-4)+2^{\frac{1-r}{2r}}}{m-2}\\
   &\qquad\leq C_{H_{r,m-1}}\bigg(\frac{p_{0}+p_{2}}{2p_{1}},\dots,
  \frac{p_{i-1}+p_{i+1}}{2p_i},\dots,\frac{p_{m-2}+p_{m}}{2p_{m-1}}\bigg)\\
   &\qquad\leq C_{H_{\infty,m-1}}\bigg(\frac{p_{0}+p_{2}}{2p_{1}},\dots,
  \frac{p_{i-1}+p_{i+1}}{2p_i},\dots,\frac{p_{m-2}+p_{m}}{2p_{m-1}}\bigg)
  =\max_{1\leq i\leq m-1}\frac{p_{i-1}+p_{i+1}}{2p_i}.
}
By taking the limit $r\to\infty$, it follows that
\Eq{*}{
   \dfrac{m-4+\sqrt{2}}{m-2}
  \leq \max_{1\leq i\leq m-1}\frac{p_{i-1}+p_{i+1}}{2p_i}.
}
In particular, with $p_i:=\sin\big(\frac{i}{m}\pi\big)$, we get that
\Eq{*}{
   \dfrac{m-4+\sqrt{2}}{m-2}\leq \cos\Big(\frac{\pi}{m}\Big),
}
which was to be shown. 

For $m=4$, both sides of the inequality are equal to $\frac{\sqrt{2}}{2}$ and hence equality holds in \eq{cos}.
\end{proof}

\section{An application of $q$-concave sequences}

In this section, we consider a selfmap of the space $\R^n$ which originates from the investigation of approximately convex real functions. Our main aim here is to prove that it has a unique fixed point.

In what follows, we will adopt the following convention: For an arbitrary sequence $a\in\S(1|n)$, let $a$ be extended to be in $\S(0|n+1)$ by setting $a_0:=0$ and $a_{n+1}:=0$. 
For $n\in\N$ and for a vector $\gamma=\big(\gamma_1,\dots,\gamma_{\lfloor\frac{n+1}{2}\rfloor}\big)\in\R^{\lfloor\frac{n+1}{2}\rfloor}$, we define the map $\T_\gamma:\R^n\to\R^n$ by
\Eq{*}{
  \big(\T_\gamma(a)\big)_i
  :=\min_{1\leq j\leq\min(i,n+1-i)}\Big(\frac{a_{i-j}+a_{i+j}}{2}+\gamma_j\Big)
  \qquad (a\in\R^n,\,i\in\{1,\dots,n\}).
}

In order to make the map $\T_\gamma$ a contraction with respect to a suitable norm on $\R^n$, we construct new norms in terms of positive sequences. Let $|\cdot|_\infty$ denote the maximum norm on $\R^n$, which is defined as $|a|_\infty:=\max_{1\leq i\leq n}|a_i|$. If $p\in\S(1|n)$ is a sequence with positive members, then we define $\|\cdot\|_p:\R^n\to\R$ by
\Eq{*}{
  \|a\|_p:=\max_{1\leq i\leq n}p_i^{-1}|a_i|=|p^{-1}a|_\infty\qquad(a\in\R^n).
}
It is easy to check that $\|\cdot\|_p$ is a norm, and hence $\R^n$ is a Banach space with respect to $\|\cdot\|_p$.

\Thm{1}{Let $p\in\S(1|n)$ be a sequence with positive members and define
\Eq{q}{
  q:=\max_{1\leq i\leq n}\frac{p_{i-1}+p_{i+1}}{2p_i}
  \qquad\mbox{and}\qquad
  q^*:=\begin{cases}
       q& \mbox{if } q\leq 1,\\
       T_{\lfloor\frac{n+1}{2}\rfloor}(q) & \mbox{if } q>1.
       \end{cases}
}
Then, for all $\gamma\in\R^{\lfloor\frac{n+1}{2}\rfloor}$, the mapping $\T_\gamma$ is $q^*$-Lipschitzian on the normed space $(\R^n,\|\cdot\|_p)$. In particular, if $p$ is strictly concave, then $\T_\gamma$ is a contraction on the normed space $(\R^n,\|\cdot\|_p)$.
}

\begin{proof}
First of all, for all $k\in\N$, we prove that the function $\min:\R^k\to\R$ is Lipschitzian with respect to the maximum norm $|\cdot|_\infty$ with Lipschitz modulus $L=1$. Indeed, if $x,y\in\R^k$, then
\Eq{*}{
   \min(x)=\min_{1\leq i\leq k}x_i
   &\leq \min_{1\leq i\leq k}\big(y_i+|x_i-y_i|\big)
   \leq \min_{1\leq i\leq k}\big(y_i+|x-y|_\infty\big)\\
   &\leq \min_{1\leq i\leq k}y_i+|x-y|_\infty=\min(y)+|x-y|_\infty.
}
Interchanging the roles of $x$ and $y$ in the above argument and then combining the two inequalities so obtained, we get that
\Eq{*}{
   \big|\min(x)-\min(y)\big|\leq |x-y|_\infty,
}
which proves our statement.

The definition of the number $q$ in \eq{q} implies that  $p\in\S(0,n+1)$ is a $q$-concave sequence, and according to \prp{Ch}, $q\geq\cos\big(\frac{\pi}{n+1}\big)$. Then, $q>\cos\big(\frac{\pi}{j}\big)$ for all $j\in\{1,\dots,n\}$. Therefore, applying the last inequality of \prp{qJC}, we obtain that, for all $i\in\{1,\dots,n\}$ and $j\in\{1,\dots,\min(i,n+1-i)\}$, \eq{ijT+} holds. Hence, on the same domain,
\Eq{Tj}{
  \frac{p_{i-j}+p_{i+j}}{2p_i}\leq T_j(q).
}
If $i\in\{1,\dots,n\}$, then $\min(i,n+1-i)\leq\frac{i+(n+1-i)}{2}=\frac{n+1}{2}$, which shows that the maximal value of $j$ is $\big\lfloor\frac{n+1}{2}\big\rfloor$. 
Therefore, \eq{Tj} implies that, for all $i\in\{1,\dots,n\}$ and $j\in\{1,\dots,\min(i,n+1-i)\}$,
\Eq{Tj+}{
  \frac{p_{i-j}+p_{i+j}}{2p_i}
  \leq \max\big\{T_1(q),\dots,T_{\lfloor\frac{n+1}{2}\rfloor}(q)\big\}.
}
In what follows, we show that the right hand side of this inequality equals $q^*$.

If $\cos\big(\frac{\pi}{n+1}\big)\leq q<1$, then $0< \arccos(q)\leq \frac{\pi}{n+1}$. Therefore, according to the first part of \lem{MT}, the sequence $T_j(q)$ is decreasing for $j\in\{0,\dots,n+1\}$ and hence $T_j(q)\leq T_1(q)=q=q^*$ for all $j\in\{1,\dots,\big\lfloor\frac{n+1}{2}\big\rfloor\}$. If $q=1$, then $T_j(q)=1=q^*$ for all $j\in\N$. On the other hand, $1<q$, then according to the second part of \lem{MT}, the sequence $(T_i(q))_{i=1}^\infty$ is increasing and hence $T_j(q)\leq T_{\lfloor\frac{n+1}{2}\rfloor}(q)=q^*$ for all $j\in\{1,\dots,\big\lfloor\frac{n+1}{2}\big\rfloor\}$.

Observe that, by the definition of the norm $\|\cdot\|_p$, for every $a\in\R^n$, we have that $|a_i|\leq p_i\|a\|_p$ is valid for $i\in\{0,1,\dots,n,n+1\}$. Now let $i\in\{1,\dots,n\}$ be fixed. Using the above established Lipschitz property of the minimum function with $k:=\min(i,n+1-i)$ and the inequality \eq{Tj+}, for all $a,b\in\R^n$, we get
\Eq{*}{
  p_i^{-1}\big|\big(\T_\gamma(a)\big)_i&-\big(\T_\gamma(b)\big)_i\big|\\
  &=p_i^{-1}\bigg|\min_{1\leq j\leq\min(i,n+1-i)}\Big(\frac{a_{i-j}+a_{i+j}}{2}+\gamma_j\Big)-\min_{1\leq j\leq\min(i,n+1-i)}\Big(\frac{b_{i-j}+b_{i+j}}{2}+\gamma_j\Big)\bigg|\\
  &\leq p_i^{-1}\max_{1\leq j\leq\min(i,n+1-i)}
  \bigg|\Big(\frac{a_{i-j}+a_{i+j}}{2}+\gamma_j\Big)-\Big(\frac{b_{i-j}+b_{i+j}}{2}+\gamma_j\Big)\bigg|\\
  &\leq\max_{1\leq j\leq\min(i,n+1-i)}
  \frac{|a_{i-j}-b_{i-j}|+|a_{i+j}-b_{i+j}|}{2p_i}\\
  &\leq \max_{1\leq j\leq\min(i,n+1-i)}
  \frac{p_{i-j}+p_{i+j}}{2p_i}\|a-b\|_p\\
  &\leq\max_{1\leq j\leq\min(i,n+1-i)}T_j(q)\|a-b\|_p
  \leq q^* \|a-b\|_p.
}
Now, upon taking the maximum with respect to $i\in\{1,\dots,n\}$, we arrive at
\Eq{*}{
  \big\|\T_\gamma(a)-\T_\gamma(b)\big\|_p\leq q^*\|a-b\|_p,
}
which completes the proof of the $q^*$-Lipschitz property of $\T_\gamma$ on $(\R^n,\|\cdot\|_p)$.

If if the sequence $p$ is strictly concave, then it is $q$-concave with some $q<1$. Therefore, the $q$-Lipschitz property of $\T_\gamma$ shows that $\T_\gamma$ is a $q$-contraction on $(\R^n,\|\cdot\|_p)$. 
\end{proof}

\Cor{Main}{For all $\gamma\in\R^{\lfloor\frac{n+1}{2}\rfloor}$, the mapping $\T_\gamma:\R^n\to\R^n$ has a unique fixed point in $\R^n$.}

\begin{proof} Let $p_i:=i(n+1-i)$ for $i\in\{0,\dots,n+1\}$. Then, by the geometric mean-arithmetic mean inequality, we have that $p_i\leq \big(\frac{n+1}{2}\big)^2$. Thus, for all $i\in\{1,\dots,n\}$ and $j\in\{1,\dots,\min(i,n+1-i)\}$, we have
\Eq{*}{
  \frac{p_{i-j}+p_{i+j}}{2p_i}
  &=\frac{(i-j)(n+1-i+j)+(i+j)(n+1-i-j)}{2i(n+1-i)}\\
  &=\frac{2i(n+1)-2i^2-2j^2}{2i(n+1-i)}
  =\frac{i(n+1)-i^2-j^2}{i(n+1-i)}
  \leq \frac{i(n+1-i)-1}{i(n+1-i)}\\
  &=\frac{p_i-1}{p_i}
  \leq \frac{\big(\frac{n+1}{2}\big)^2-1}{\big(\frac{n+1}{2}\big)^2}
  =\frac{n^2+2n-3}{n^2+2n+1}
  =\frac{(n-1)(n+3)}{(n+1)^2}.
}
Therefore, the sequence $p\in\S(0,n+1)$ is $q$-concave with $q=\frac{(n-1)(n+3)}{(n+1)^2}<1$. According to the \thm{1}, the mapping $\T_\gamma$ is a $q$-contraction on $(\R^n,\|\cdot\|_p)$. Therefore, by the Banach Fixed Point theorem, it possesses a unique fixed point.
\end{proof}


\begin{thebibliography}{10}

\bibitem{HarLitPol34}
G.~H. Hardy, J.~E. Littlewood, and G.~Pólya, \emph{{Inequalities}}, Cambridge
  University Press, Cambridge, 1934, (first edition), 1952 (second edition).
  \MR{13,727e}

\bibitem{Kra16}
X.~Z. Krasniqi, \emph{{On {$\alpha$}-convex sequences of higher order}}, J.
  Numer. Anal. Approx. Theory \textbf{45} (2016), no.~2, 177–182.
  \MR{3599323}

\bibitem{Kuc85}
M.~Kuczma, \emph{{An {I}ntroduction to the {T}heory of {F}unctional {E}quations
  and {I}nequalities}}, {Prace Naukowe Uniwersytetu Śląskiego w Katowicach},
  vol. 489, Państwowe Wydawnictwo Naukowe — Uniwersytet Śląski,
  Warszawa–Kraków–Katowice, 1985, 2nd edn. (ed. by A. Gilányi),
  Birkhäuser, Basel, 2009. \MR{0788497 (86i:39008), MR 2467621}

\bibitem{Mit70}
D.~S. Mitrinović, \emph{{Analytic inequalities}}, {Die Grundlehren der
  mathematischen Wissenschaften, Band 165}, Springer-Verlag, New York-Berlin,
  1970, In cooperation with P. M. Vasić. \MR{0274686}

\bibitem{MitPecFin91}
D.~S. Mitrinović, J.~E. Pečarić, and A.~M. Fink, \emph{{Inequalities
  {I}nvolving {F}unctions and {T}heir {I}ntegrals and {D}erivatives}},
  {Mathematics and its Applications (East European Series)}, vol.~53, Kluwer
  Academic Publishers Group, Dordrecht, 1991. \MR{93m:26036}

\bibitem{MitPecFin93}
D.~S. Mitrinović, J.~E. Pečarić, and A.~M. Fink, \emph{{Classical and {N}ew {I}nequalities in {A}nalysis}},
  {Mathematics and its Applications (East European Series)}, vol.~61, Kluwer
  Academic Publishers Group, Dordrecht, 1993. \MR{94c:00004}

\bibitem{NicPer06}
C.~P. Niculescu and L.-E. Persson, \emph{{Convex {F}unctions and {T}heir
  {A}pplications}}, {CMS Books in Mathematics/Ouvrages de Mathématiques de la
  SMC, 23}, Springer-Verlag, New York, 2006, A contemporary approach.

\bibitem{Nie11}
M.~Niezgoda, \emph{{Remarks on convex functions and separable sequences,
  {II}}}, Discrete Math. \textbf{311} (2011), no.~2-3, 178–185. \MR{2739922}

\bibitem{Nie17b}
M.~Niezgoda, \emph{{Inequalities for convex sequences and nondecreasing convex functions}}, Aequationes Math. \textbf{91} (2017), no.~1, 1–20.
  \MR{3600782}

\bibitem{Nie17a}
M.~Niezgoda, \emph{{Sherman, {H}ermite-{H}adamard and {F}ejér like inequalities for convex sequences and nondecreasing convex functions}}, Filomat \textbf{31} (2017), no.~8, 2321–2335. \MR{3637029}

\bibitem{Pop44}
T.~Popoviciu, \emph{{Les fonctions convexes}}, Hermann et Cie, Paris, 1944.
  \MR{8,319a}

\bibitem{RobVar73}
A.~W. Roberts and D.~E. Varberg, \emph{{Convex {F}unctions}}, {Pure and Applied Mathematics}, vol.~57, Academic Press, New York–London, 1973. \MR{56
  \#1201}

\bibitem{SofTinAcu18}
D.~F. Sofonea, I.~Ţincu, and A.~M. Acu, \emph{{Convex sequences of higher
  order}}, Filomat \textbf{32} (2018), no.~13, 4655–4663. \MR{3897332}

\bibitem{TabTabZol12}
Ja. Tabor, Jó. Tabor, and M.~Żołdak, \emph{{Strongly convex sequences}},
  {Inequalities and applications 2010}, {Internat. Ser. Numer. Math.}, vol.
  161, Birkhäuser/Springer, Basel, 2012, p.~183–188. \MR{3203786}

\bibitem{WuDeb07}
Sh. Wu and L.~Debnath, \emph{{Inequalities for convex sequences and their
  applications}}, Comput. Math. Appl. \textbf{54} (2007), no.~4, 525–534.
  \MR{2340843}

\bibitem{Yil18}
Ş. Yıldız, \emph{{A general matrix application of convex sequences to
  {F}ourier series}}, Filomat \textbf{32} (2018), no.~7, 2443–2449.
  \MR{3900949}

\end{thebibliography}

\providecommand{\bysame}{\leavevmode\hbox to3em{\hrulefill}\thinspace}
\providecommand{\MR}{\relax\ifhmode\unskip\space\fi MR }
\providecommand{\MRhref}[2]{%
  \href{http://www.ams.org/mathscinet-getitem?mr=#1}{#2}
}
\providecommand{\href}[2]{#2}

\end{document}